\documentclass[12pt,a4paper]{amsart}

\usepackage{enumerate}
\usepackage{verbatim}
\usepackage[all]{xy}

\usepackage[pdftex,pagebackref]{hyperref}
\hypersetup{letterpaper=true,colorlinks=true}
\hypersetup{pdfpagemode=UseNone,urlcolor=blue}
\hypersetup{linkcolor=blue,citecolor=blue,pdfstartview=FitH}

\usepackage{amsmath,amsfonts,amssymb,amsthm}
\usepackage{graphicx}

\usepackage{color}

\usepackage{ulem}

\def\emph#1{\textit{#1}}

\renewcommand{\geq}{\geqslant}
\renewcommand{\leq}{\leqslant}

\renewcommand{\bar}{\overline}

\def\Z{{\mathbb Z}}

\newtheorem{theorem}{Theorem}[section]
\newtheorem{lemma}[theorem]{Lemma}

\newtheorem{corollary}[theorem]{Corollary}
\newtheorem{proposition}[theorem]{Proposition}

\newtheorem{conjecture}[theorem]{Conjecture}

\renewcommand{\hom}{\textrm{Hom}}

\newcommand{\pic}[1]{\textrm{Pic}(#1)}   
\newcommand{\slope}[1]{\frac{\deg #1}{\rk #1}}   
\newcommand{\coker}{\textrm{coker}}   

\newcommand{\ext}{\textrm{Ext}}

\newcommand{\ch}[1]{\textrm{ch}(#1)}   
\newcommand{\td}[1]{\textrm{td}(#1)}   
\newcommand{\rk}{\textrm{rank}}
\newcommand{\im}{\textrm{Im}}

\newcommand{\fM}{\mathfrak{M}}   

\newcommand{\cX}{\mathcal{X}}

\newcommand{\cL}{\mathcal{L}}
\newcommand{\cext}{\mathcal{E}\mathit{xt}}
\newcommand{\chom}{\mathcal{H}\mathit{om}}

\newenvironment{acknowledgements}{\textbf{Acknowledgements:}}{}

\title{A Note on Singular Moduli Spaces of Sheaves on K3 Surfaces}
\date{\today}
\author{Ziyu Zhang}
\address{Institute for Mathematics, University of Mainz, Staudingerweg 9, 55099 Mainz, Germany}
\email{zhangzy@uni-mainz.de}

\begin{document}

\maketitle

\begin{abstract}
This paper studies deformations and birational maps between singular moduli spaces of semistable sheaves with 2-divisible Mukai vectors on K3 surfaces. It is showed that under certain conditions, two such moduli spaces of the same dimension can be connected by deformations and birational maps.
\end{abstract}

\tableofcontents

\section{Introduction}

Moduli spaces of semistable sheaves have been studied for a long time. Let $X$ be a K3 surface, and $H\in \pic{X}$ be an ample line bundle. Let $v\in H^{\textrm{even}}(X,\Z)$ be a fixed Mukai vector. Then there is a moduli space $\fM_{X,H}(v)$ which parameterizes S-equivalence classes of semistable sheaves with respect to the polorization $H$, whose Mukai vectors are $v$. These moduli spaces were first constructed by Gieseker \cite{Gies} and Maruyama \cite{Maru1,Maru2}, and then studied by many other people. 

When the polarization $H$ is generic and the Mukai vector $v$ is primitive, namely, the greatest common divisor of all components of $v$ is 1, every semistable sheaf must be stable. In this case, Mukai \cite[Corollary 0.2]{Mukai-symp} proved that the muduli space $\fM_{X,H}(v)$ is a smooth irreducible holomorphic symplectic manifold. 

It is an interesting problem to study the relation among these smooth moduli spaces. There is a very nice theorem on irreducible holomorphic symplectic manifolds due to Huybrechts \cite[Theorem 2.5]{Huy-Kahlercone}, which states that two birational irreducible holomorphic symplectic manifolds are always deformation equivalent. Applying this theorem, Yoshioka proved that \cite[Theorem 8.1]{Yosh-abelian}, if two such smooth moduli spaces have the same dimension, then they are deformation equivalent. This is really a nice result. However, it's not good news in the study of irreducible holomorphic symplectic manifolds. The reason is, although we have lots of choices for the underlying K3 surface, the generic polarization and the primitive Mukai vector, the resulting moduli spaces provide only one deformation type of holomorphic symplectic manifolds in every even dimension. 

In this paper, we are trying to generalize Yoshioka's result to the case of 2-divisible Mukai vectors, which are Mukai vectors whose greatest common divisor among all components is 2. More precisely, we will prove the following theorem: 

\begin{theorem}\label{_thm_main_}
For $i=1$ or $2$, assume $X_i$ is a projective K3 surface,
$v_i=(r_i,c_i,a_i)$ is a primitive Mukai vector with $r_i>0$. $H_i$ is a polarization on $X_i$ which is generic
with respect to the Mukai vector $2v_i$. Assume further that $r_i$ and
$c_i$ are coprime. If $\dim \fM_{X_1,H_1}(2v_1)=\dim \fM_{X_2,H_2}(2v_2)$,
then $\fM_{X_1,H_1}(2v_1)$ and $\fM_{X_2,H_2}(2v_2)$ can be connected by a series of
deformations and birational maps.
\end{theorem}

The proof contains two main steps. In the first step, we will show that, for two such moduli spaces $\fM_{X_1,H_1}(2v_1)$ and $\fM_{X_2,H_2}(2v_2)$, if they have the same dimension and $r_1=r_2$, then they are deformation equivalent. The idea of this part of proof is to use the deformations of polarized K3 surfaces to deform both of the given moduli spaces to a third moduli space of sheaves on an elliptic K3 surface with a rank 2 Picard lattice. So that we know all moduli spaces in the same dimension with the same rank parameter are deformation equivalent. In the second step, we fix a dimension of the moduli spaces and let rank vary. For every possible value of the rank component of the Mukai vector which is at least 8, we find one particular moduli space of sheaves on a K3 surface with rank 1 Picard lattice, and prove it is birational to a certain moduli space of sheaves of rank 2. The birational map between these two moduli spaces is established via extensions by exception bundles, which was introduced in \cite{Yosh-exceptional}. Combining the two steps, any two moduli spaces of sheaves of rank not equal to 4 or 6, as stated in theorem \ref{_thm_main_}, can be connected in at most three steps, namely a deformation followed by a birational map, then by another deformation. For a technical reason, we have to deal with the case of rank 4 and 6 seperately. To connect these remaining moduli spaces, we just need to realize that, in every dimension we can find a moduli space of sheaves of rank 4 (respectively 6), which is birational to another moduli space of rank 14 (respectively 20). So now all moduli spaces of the same dimension are connected by deformations and birational maps. 

However, due to the lack of a version of Huybrechts's theorem \cite[Theorem 2.5]{Huy-Kahlercone} in the singular case, here we cannot get rid of the birational maps and get a theorem as nice as Yoshioka's \cite[Theorem 8.1]{Yosh-abelian}. However, we hope a similar result of deformation equivalence of singular moduli spaces in the same dimension is still true for arbitrary non-primitive Mukai vectors. So we can formulate the following conjecture:

\begin{conjecture}\label{def-conj}
Let $m$ be a positive integer. For $i=1$ or $2$, assume $X_i$ is a
projective K3 surface, $v_i=(r_i,c_i,a_i)$ is a primitive positive
Mukai vector, and $H_i$ is a polarization which is generic with
respect to the Mukai vector $mv_i$. If $\dim M_{X_1,H_1}(mv_1)=\dim
M_{X_2,H_2}(mv_2)$, then $M_{X_1,H_1}(mv_1)$ and $M_{X_2,H_2}(mv_2)$
are deformation equivalent.
\end{conjecture}

This paper is organized as follows: the two steps in the proof of theorem \ref{_thm_main_} will be done in sections 2 and 3. In the proof of the deformation equivalence we will need a result of local finiteness of walls in the ample cone (lemma \ref{2times}), whose proof will be given in the end of section 2. In section 3, for every $r\geq 4$ we will find a moduli space of sheaves of rank $2r$, and prove it's birational to a moduli space of sheaves of rank 2. The two remaining cases for $r=2$ and $r=3$ will be treated at the end of the section (propositions \ref{exp1} and \ref{exp2}). We will omit the proofs because they will be identical with that of the proposition \ref{rk3}. These two sections finish the proof of Theorem \ref{_thm_main_}. Finally, in section 4, we will prove that certain moduli spaces of torsion sheaves with 2-divisible Mukai vectors are birational to moduli spaces of sheaves of positive ranks. 

\begin{acknowledgements}
I would like to express my deep gratitude to my PhD advisor Jun Li, for all his support and encouragement. I would also like to thank K\={o}ta Yoshioka for kindly answering my questions and pointing out references, and thank Jason Lo for helpful discussions. I also appreciate the help of Daniel Huybrechts and Manfred Lehn in the final stage of this work. 
\end{acknowledgements}

\section{Deformation Equivalence}

In this section, we will show that any two moduli spaces of sheaves on K3 surfaces with the same dimension, rank and divisibility are deformation equivalent to each other. This will be obtained by deforming the two moduli space to a third moduli space of sheaves of the same rank over an elliptic K3 surface. More precisely, we will show: 

\begin{proposition}\label{deform}
Let $X$ be a K3 surface, $v=(r, c, a)\in H^{even}(X)$ be a Mukai vector with $\gcd (r, c)=1$, $H$ be a generic polarization. Let $X^e$ be an elliptic K3 surface with $\pic{X^e}=\mathbb{Z}[\sigma]\oplus\mathbb{Z}[f]$ where $\sigma$ is the class of a section of the elliptic fibration and $f$ is the fiber class. The the moduli space $\fM_{X,H}(2v)$ is deformation equivalent to a certain moduli space of sheaves $\fM_{X^e,\sigma+lf}(2v')$ on $X^e$, where $v'=(r, \sigma+lf, a')$ for some $l$ and $a'$.
\end{proposition}

Before proving this proposition, we have to consider moduli spaces of polarized K3 surfaces. Much of this has been summarized in \cite{HuyLehn}.

Let $d$ be a positive number, we consider all K3 surfaces $X$ with ample primitive line bundles $L$ satisfying $c_1^2(L)=2d$. Then there is a quasi-projective scheme $\mathcal{K}_d$, which is a coarse moduli space of all such pairs $(X, L)$. 

This moduli space can be constructed by GIT. More precisely, it is a $\textrm{PGL}(N)$ quotient of $\mathcal{H}_d$ which is an open subset of a certain Hilbert scheme. The universal family over the Hilbert scheme provides a universal family of polarized K3 surfaces over $\mathcal{H}_d$. For every point $t\in\mathcal{H}_d$, the fiber of this universal family $(\cX _t, \cL _t)$ is exactly the pair corresponding to the image of $t$ in $\mathcal{K}_d$. 

Furthermore, we know that both $\mathcal{K}_d$ and $\mathcal{H}_d$ are irreducible. Therefore, any two primitively polarized K3 surfaces $(X_1, L_1)$ and $(X_2, L_2)$ with $L_1^2=L_2^2$ are deformation equivalent to each other. 

Another fact which will be used later is: for a general polarized K3 surface $(X, L)\in\mathcal{K}_d$, we have $\pic{X}=\mathbb{Z}[H]$. In other words, in the moduli space $\mathcal{K}_d$ of primitively polarized K3 surfaces, away from a countable union of Zariski closed subsets, the Picard number $\rho(X)=1$. However, the countable union of polarized K3 surfaces $(X, H)\in\mathcal{K}_d$ with $\rho(X)\geq 2$ is also dense in $\mathcal{K}_d$. 

We will need the following proposition: 

\begin{proposition}\label{rel-moduli}
Let $r$ and $k$ be positive integers and $a$ be an arbitrary integer. Then there is a relative moduli space $\varphi: \fM\longrightarrow \mathcal{H}_d$ of semistable sheaves, such that for every $t\in\mathcal{H}$, $\varphi^{-1}(t)$ is isomorphic to the moduli space of semistable sheaves $\fM_{\cX_t, \cL_t}(2v_t)$ where $v_t=(r, k\cL_t, a)$. 
\end{proposition}

\begin{proof}
The proposition is a special case of the existence theorem of relative moduli spaces of semistable sheaves (\cite[Theorem 4.3.7]{HuyLehn}).
\end{proof}

Now we are ready to prove proposition \ref{deform}. The proof will be parallel to the proof of Theorem 6.2.5 in \cite{HuyLehn}. 

\begin{proof}[Proof of Proposition \ref{deform}]
First of all, without loss of generality, we can assume that $c$ is an ample class. 

In fact, if $c$ is not ample, we can always twist the sheaves by the ample line bundle $H$ sufficiently many times to get a new moduli space isomorphic to the original one. More precisely, for every semistable sheaf $F$ with Mukai vector $2v=2(r, c, a)$, we consider a new sheaf $F\otimes H^{\otimes m}$. We have that $$v(F\otimes H^{\otimes m})=v(F)\cdot\mathrm{ch}(H)^{\otimes m}=2(r, c+rmH, a+cmH+\frac{rm^2}{2}H^2).$$ We denote the new Mukai vector $(r, c+rmH, a+cmH+\frac{rm^2}{2}H^2)$ by $v'$. Note that the above tensoring procedure can also be done in families, and is invertible by tensoring negative powers of the line bundle $H$. Therefore we have an isomorphism between $\fM_{X,H}(2v)$ and $\fM_{X,H}(2v')$. Since $H$ is ample, we know that $c+rmH$ is also an ample class when $m$ is sufficiently large. Therefore we can replace $\fM_{X,H}(2v)$ by $\fM_{X,H}(2v')$. 

Secondly, without loss of generality, we can assume that $\rho(X)\geq 2$. 

If not, then $\pic{X}=\mathbb{Z}[H]$. Let $c=kH$ where $k>0$. We can consider the relative moduli space $\varphi:\fM\longrightarrow\mathcal{H}_d$ in proposition \ref{rel-moduli}. Our moduli space $\fM_{X,H}(2r, 2kH, 2a)$ is a fiber of this relative moduli space over the point $(X, H)$. Since polarized K3 surfaces with Picard number at least two are dense in $\mathcal{H}_d$, a priori, we can choose any polarized K3 surface $(X', H')$ near $(X, H)$ in $\mathcal{H}_d$, with $\rho(X')\geq 2$, so that the original moduli space $\fM_{X,H}(2r, 2kH, 2a)$ is deformation equivalent to the new moduli space $\fM_{X',H'}(2r, 2kH', 2a)$. 

However, we still have to make sure that $H'$ is a generic polarization on $X$ with respect to the Mukai vector $(2r, 2kH', 2a)$. We need the following lemma: 

\begin{lemma}\label{2times}
In a sufficiently small open neighborhood of $(X, H)$ in $\mathcal{H}_d$, there are at most finitely many hypersurfaces, such that for every point $(X', H')$ not on those hypersurfaces, the polarization $H'$ is generic with respect to the corresponding Mukai vector $2(r, kH', a)$. 
\end{lemma}

We will postpone the proof of this lemma to the end of this section. 

By virtual of the above lemma, we can always find a polarized K3 surface $(X', H')$ with $\rho(X')\geq 2$, which is deformation equivalent to the original polarized K3 surface $(X, H)$. 

Thirdly, we can even assume that $c$ and $H$ are linearly dependent, and $H^2\gg 0$. 

In fact, since $H$ is in an open chamber and $\rho(X)\geq 2$, we can always pick another primitive ample line bundle $H'$ in the same open chamber as $H$, which is linearly independent with $c$. We can now twist the sheaves in the moduli space $\fM_{X,H}(2r, 2c, 2a)$ by $H'$ for a few times to get an isomorphic moduli space. More precisely, if $v(F)=2(r, c, a)$, then $$v(F\otimes H'^{\otimes m})=v(F)\cdot\mathrm{ch}(H')^{\otimes m}=2(r, c+rmH', a+cmH'+\frac{rm^2}{2}H'^2).$$ We denote the Mukai vector $(r, c+rmH', a+cmH'+\frac{rm^2}{2}H'^2)$ by $v'$. So $\fM_{X,H}(2v)$ is isomorphic to $\fM_{X,H}(2v')$. Note that when $m$ is sufficiently large, $c+rmH'$ is very close to $H'$, therefore is also ample and in the same chamber as $H'$, or $H$. So the moduli space $\fM_{X,H}(2v')$ is isomorphic to $\fM_{X,c+rmH'}(2v')$. With the additional assumption that $r$ and $c$ are coprime, we know that $c+rmH'$ is primitive. Increase $m$ we can make $(c+rmH')^2\gg 0$. 

Finally, we can finish the proof based on the above three steps of reductions. 

Note that there is an elliptic K3 surface $X^e$ as stated in the proposition. On $X^e$, the class $\sigma +lf$ is always ample and suitable when $l$ is sufficiently large. We choose the value of $l$ by the equation $H^2=(\sigma +lf)^2$. Then the polarized K3 surface $(X^e, \sigma +lf)$ is in $\mathcal{H}_d$. By proposition \ref{rel-moduli}, we know that $\fM_{X,H}(2(r, c, a))$ is deformation equivalent to $\fM_{X^e,\sigma +lf}(2(r, \sigma +lf, a))$. 
\end{proof}

From proposition \ref{deform} we can prove the following corollary:

\begin{corollary}\label{_def_cor_}
For $i=1,2$, let $X_i$ be a K3 surface, $v_i=(r_i, c_i, a_i)\in H^{\textrm{even}}(X_i, \mathbb{Z})$ be a primitive Mukai vector with $\mathrm{gcd}(r_i, c_i)=1$ and $H_i$ be a generic ample line bundle on $X_i$ with respect to $v_i$. If $r_1=r_2=r$ and $\langle v_1^2\rangle=\langle v_2^2 \rangle$, then moduli spaces of semistable sheaves $\fM_{X_1, H_1}(2v_1)$ is deformation equivalent to $\fM_{X_2, H_2}(2v_2)$. 
\end{corollary}

\begin{proof}
By proposition \ref{deform}, there exists an elliptic K3 surface $X^e$, such that, for $i=1$ and $2$, there exists an positive integer $l_i$, such that $\fM_{X_i, H_i}(2v_i)$ is deformation equivalent to $\fM_{X^e, \sigma +l_if}(2(r, \sigma +l_if, a_i))$. 

Now we claim the two moduli spaces $$\fM_{X^e, \sigma +l_1f}(2(r, \sigma +l_1f, a_1))\quad\mbox{and}\quad\fM_{X^e, \sigma +l_2f}(2(r, \sigma +l_2f, a_2))$$ are isomorphic. In fact, from the given condition we know the two moduli spaces have the same dimension, therefore $\langle (r, \sigma +l_1f, a_1)^2\rangle=\langle (r, \sigma +l_2f, a_2)^2\rangle $, that is, $(\sigma +l_1f)^2-2ra_1=(\sigma +l_2f)^2-2ra_2$. This implies $2r(a_2-a_1)=2(l_2-l_1)\sigma\cdot f=2(l_2-l_1)$, so $r$ divides $l_2-l_1$. Let $l_2-l_1=r\cdot h$. 

It's now easy to see that, there is an isomorphism from $\fM_{X^e, \sigma +l_1f}(2(r, \sigma +l_1f, a_1))$ to $\fM_{X^e, \sigma +l_2f}(2(r, \sigma +l_2f, a_2))$. For every sheaf $F\in\fM_{X^e, \sigma +l_1f}(2(r, \sigma +l_1f, a_1))$, we have $F\otimes \mathcal{O}(f)^{\otimes h}\in\fM_{X^e, \sigma +l_2f}(2(r, \sigma +l_2f, a_2))$. Although $\sigma +l_1f$ and $\sigma +l_2f$ might be different polarizations, they are both in the open chamber containing the fiber class $f$. Therefore, the two moduli spaces are isomorphic. 

This proves the two original moduli spaces are deformation equivalent to each other. 
\end{proof}

To conclude this section we prove the Lemma \ref{2times}. The idea is similar to Lemma 4.C.2 in \cite{HuyLehn}. For the proof we need to state the following lemma: 

\begin{lemma}\cite{HuyLehn}[Theorem 4.C.3]\label{polar}
Let $H$ be an ample line bundle, and $F$ be a semistable coherent sheaf with Mukai vector $v(F)=(r, c, a)$. Let $F'$ be a subsheaf of $F$ with Mukai vector $v(F')=(r', c', a')$, where $0<r'<r$, such that the reduced Hilbert polynomials $p(F, m)=p(F', m)$. Let $\Delta=2r^2+c^2-2ra$ be the discriminant of $F$. Then we have
\begin{enumerate}
\item $ra'-r'a=0$;
\item let $\xi=rc'-r'c$, then $\xi\cdot H=0$;
\item $-\frac{r^2}{4}\Delta\leq\xi^2\leq 0$ and $\xi^2=0$ if and only if $\xi=0$. \qed
\end{enumerate}
\end{lemma}

Now we can begin to prove lemma \ref{2times}. 

\begin{proof}[Proof of Lemma \ref{2times}]
We consider the walls in $H^2(X,\mathbb{R})$. Let $\sigma$ be a generator of $H^{2,0}(X)$, then $\bar{\sigma}$ is a generator of $H^{0,2}(X)$. Let $$e_1=\frac{\sigma+\bar{\sigma}}{\sqrt{(\sigma+\bar{\sigma})^2}}$$ and $$e_2=\frac{\sigma-\bar{\sigma}}{\sqrt{-(\sigma-\bar{\sigma})^2}}.$$ Then both $e_1$ and $e_2$ are unit vectors in $H^2(X,\mathbb{R})$. Furthermore, since the signature of Poincare pairing on $H^2(X,\mathbb{R})$ is $(3, h^{1,1}-3)$, we see that $\{e_1, e_2, H\}$ is an orthogonal basis of a three dimensional subspace $V$ of $H^2(X,\mathbb{R})$. Therefore, for every class $u\in H^2(X,\mathbb{R})$, we can decompose it into $u=a_1e_1+a_2e_2+a_0H+u_0$ where $a_1, a_2\in\mathbb{R}$ and $u_0\in V^{\perp}$. Besides the Poincare pairing, we can define an Euclidean pairing on $H^2(X,\mathbb{R})$, by $\parallel u\parallel^2=\sqrt{a_1^2+a_2^2+a_0^2H^2-u_0^2}$. 

For every pairing $(X', H')\in\mathcal{H}_d$, we can produce $e'_1, e'_2$ in the same way as $e_1, e_2$. In a sufficiently small open neighborhood of the pair $(X,H)$ in $\mathcal{H}_d$, every base manifold $X'$ is a small deformation of $X$, whose $H^{2,0}$ and $H^{0,2}$ are very closed to $H^{2,0}(X)$ and $H^{0,2}(X)$. More precisely, we can require $\parallel e'_1-e_1\parallel<\varepsilon$ and  $\parallel e'_2-e_2\parallel<\varepsilon$. 

Since $e'_1$ and $e'_2$ are both perpendicular to $H$, let $e'_1=\lambda_1e_1+\delta_1e_2+h_1$ where $h_1\in V^{\perp}$ and $e'_2=\delta_2e_1+\lambda_2e_2+h_2$ where $h_2\in V^{\perp}$. Then $\parallel e'_1-e_1\parallel<\varepsilon$ implies $(1-\lambda_1)^2+\delta_1^2+\parallel 	h_1\parallel^2<\varepsilon ^2$. So we can conclude $|1-\lambda_1|<\varepsilon$, $|\delta_1|<\varepsilon$, $\parallel h_1\parallel<\varepsilon$. Similarly, from $\parallel e'_1-e_1\parallel<\varepsilon$, we can deduce $|1-\lambda_2|<\varepsilon$, $|\delta_2|<\varepsilon$, $\parallel h_2\parallel<\varepsilon$. 

By lemma \ref{polar}, we know that there exist a positive constant $\Delta$, such that for every $\xi\in H^{1,1}(X')$ which produces a wall in $H^2(X', \mathbb{R})=H^2(X, \mathbb{R})$ for some pair $(X', H')$ in the small neighborhood of $(X, H)$, we have $-\Delta<\xi^2\leq 0$. Since $\xi\perp H$, we can decompose $\xi$ as $\xi=a_1e_1+a_2e_2+\xi_0$ where $\xi_0\in V^{\perp}$. We want to show the Euclidean norm $\parallel\xi\parallel^2=a_1^2+a_2^2-\xi_0^2$ is bounded. 

Note that $\xi\perp e'_1$ and $\xi\perp e'_2$, that is, 
\begin{eqnarray*}
\lambda_1a_1+\delta_1a_2+h_1\cdot\xi_0=0\\
\delta_2a_1+\lambda_2a_2+h_2\cdot\xi_0=0
\end{eqnarray*}
Therefore we have
\begin{eqnarray*}
(\lambda_1a_1+\delta_1a_2)^2=(h_1\cdot\xi_0)^2\leq\parallel h_1\parallel^2\parallel\xi_0\parallel^2\\
(\delta_2a_1+\lambda_2a_2)^2=(h_2\cdot\xi_0)^2\leq\parallel h_2\parallel^2\parallel\xi_0\parallel^2
\end{eqnarray*}
We add the two inequalities. The left hand side is 
\begin{eqnarray*}
& &(\lambda_1^2+\delta_2^2)a_1^2+(\delta_1^2+\lambda_2^2)a_2^2 +2(\lambda_1\delta_1+\delta_2\lambda_2)a_1a_2\\
&>&(1-\varepsilon)^2(a_1^2+a_2^2)-4\varepsilon(1+\varepsilon)a_1a_2\\
&>&(1-\varepsilon)^2(a_1^2+a_2^2)-2\varepsilon(1+\varepsilon)(a_1^2+a_2^2)\\
&=&(1-4\varepsilon-2\varepsilon^2)(a_1^2+a_2^2).
\end{eqnarray*}
So now we have $$(1-4\varepsilon-2\varepsilon^2)(a_1^2+a_2^2)<2\varepsilon^2\parallel\xi_0\parallel^2.$$
Therefore we have $$-\Delta<a_1^2+a_2^2-\parallel\xi_0\parallel^2 <(\frac{2\varepsilon}{1-4\varepsilon-2\varepsilon^2}-1)\parallel\xi_0\parallel^2.$$
So $$\parallel\xi_0\parallel^2< \frac{\Delta}{1-\frac{2\varepsilon}{1-4\varepsilon-2\varepsilon^2}}.$$
Therefore $$\parallel\xi\parallel^2=a_1^2+a_2^2+\parallel\xi_0\parallel^2\leq 2\parallel\xi_0\parallel^2 <\frac{2\Delta}{1-\frac{2\varepsilon}{1-4\varepsilon-2\varepsilon^2}}$$ which is bounded regardless of the underlying K3 surface $X'$. 

Therefore, in a small open neighborhood of $(X,H)$ in $\mathcal{H}_d$, there are at most finitely many hypersurfaces, on which the polarization is not generic. Since we have assumed that $H$ is in generic on $X$, we can always shrink the open neighborhood of $(X,H)$ in $\mathcal{H}_d$, so that it misses all these hypersurfaces. In other words, for every point $(X', H')$ in this open neighborhood, $H'$ is generic with respect to the Mukai vector $2(r, kH', a)$. 
\end{proof}

\section{Birational Equivalence}


In this section, we will establish the birational equivalence between some singular moduli spaces of semistable sheaves, so that together with the deformation result we obtained in the previous section, for any two moduli spaces of sheaves in the same dimension, both with 2-divisible Mukai vectors, we can connect them by a series of deformations and birational maps. For this purpose, in this section we always fix the underlying K3 surface, which is a projective K3 surface with Picard number 1. Roughly speaking, we will prove that for every value of $r$, there's a moduli space of rank $2r$ sheaves, which is birational to a moduli space of rank 2 sheaves (in case $r=2$ or $3$, birational to a moduli of sheaves of higher rank). The method we are using here in proving the birationality is mainly the technique of exceptional bundles introduced by Yoshioka in \cite{Yosh-exceptional}. Similar to \cite{Yosh-exceptional}, the proof here involves some delicate analysis of slopes of sheaves. 

The main theorem of this section is:

\begin{proposition}\label{rk3}
For any given $r\geq 4$ and $s\geq 1$, let $X$ be a projective K3 surface with $\pic{X}=\mathbb{Z}[H]$ where $H$ is an ample line bundle with $H^2=2s(r-1)^2-2r$. Consider Mukai vectors
\begin{eqnarray*}
v_0&=&(r-1, (r-2)H, (r^3-5r^2+8r-4)s+(-r^2+3r-1)),\\
v_1&=&(r, (r-1)H, (r^3-4r^2+6r-4)s+(-r^2+2r-1)),\\
v_2&=&(1, H, (r^2-2r)s-r).
\end{eqnarray*}
Then the moduli spaces of semi-stable sheaves $\fM_{X,H}(2v_1)$ is birational equivalent to $\fM_{X,H}(2v_2)$.
\end{proposition}

Simple calculation shows that $\langle v_1^2\rangle =\langle v_2^2\rangle=2s$, therefore moduli spaces $\fM_{X,H}(2v_1)$ and $\fM_{X,H}(2v_2)$ have the same dimension. We also have $\langle v_0, v_1\rangle=-1$ and $\langle v_0, v_2\rangle=1$. Furthermore, we have $\langle v_0^2\rangle=-2$. Due to the following lemma of Yoshioka, the technique of exceptional bundles applies. 

\begin{lemma}\cite[Theorem 3.6]{Yosh-exceptional}
The moduli space $\fM_{X,H}(v_0)$ consists of one point which corresponds to a $\mu$-stable locally free sheaf $E$. We call such a sheaf an exceptional sheaf. \qed
\end{lemma}

We should also note that, in proposition \ref{rk3}, we require $r\geq 4$. The reason is that, one step in the proof of lemma \ref{3rdlemma} fails to go through when $r=2$ or $3$. However, by carefully choosing the numerical data in these two cases, the proof will work in exactly the same way. We will state a theorem in the two exceptional cases in the end of this section and omit its proof. 


Before we prove proposition \ref{rk3}, we want to show that a generic point in the moduli spaces $\fM_{X,H}(2v_1)$ or $\fM_{X,H}(2v_2)$ is represented by a $\mu$-stable locally free sheaf. We recall the following lemma due to Yoshioka: 

\begin{lemma}\cite[Remark 2.2]{Yosh-irr}\cite[Lemma 4.4, Remark 4.3]{Yosh-exceptional}
Let $(X,H)$ be a generic polarized K3 surface. Let $v=(lr, l\xi, a)$ be a Mukai vector, such that $r$ and $\xi$ are coprime. Then there's at least one $\mu$-stable sheaf with Mukai vector $v$ unless the following two conditions simultaneously hold: 
\begin{itemize}
\item $\displaystyle \frac{\xi^2+2}{2r}$ is an integer;
\item $\langle v^2 \rangle <2l^2$. \qed
\end{itemize}
\end{lemma}

We apply the above lemma in our situation. Note that when $v=2v_1$ or $2v_2$, we always have $$\langle v^2 \rangle=8s\geq 8=2l^2.$$ The second condition in the above lemma fails, both of the moduli spaces $\fM_{X,H}(2v_1)$ and $\fM_{X,H}(2v_2)$ contain at least one $\mu$-stable sheaf. Take into consideration that the $\mu$-stability is an open condition, we know that each of the two moduli spaces has an open subscheme, which parametrizes $\mu$-stable sheaves. Moreover, the following lemma due to Yoshioka shows that a generic point in the $\mu$-stable locus is represented by a locally free sheaf:

\begin{lemma}\cite[Remark 3.2]{Yosh-abelian}\label{st+lf}
Let $r$ be the rank component of the Mukai vector $v$. Then the complement of the locus of locally free sheaves in the moduli of $\mu$-stable sheaves $\fM_{X,H}^{\mu\textrm{-st}}(v)$ has codimension $r-1$. \qed
\end{lemma}

In our situation, the rank of the sheaf is at least 2. So we conclude that

\begin{corollary}
There are open subschemes $U_1$ and $U_2$ in $\fM_{X,H}(2v_1)$ and $\fM_{X,H}(2v_2)$ respectively, which parametrize all locally free $\mu$-stable sheaves with Mukai vector $2v_1$ and $2v_2$. \qed
\end{corollary}


First of all we need the following lemma, whose proof is very similar to \cite[Lemma 2.1]{Yosh-exceptional}.

\begin{lemma}\label{l3-1}
For any $\mu$-stable locally free sheaf $G$ in $\fM_H(2v_2)$ and any non-trivial extension $$0\longrightarrow E\longrightarrow V\longrightarrow G\longrightarrow 0,$$ $V$ is also a $\mu$-stable locally free sheaf. 
\end{lemma}
\begin{proof}
Since both $E$ and $G$ are locally free, $V$ is also locally free. Assume $V$ is not $\mu$-stable, then we can find a locally free subsheaf $K$ such that $$\slope{K}\geq\slope{V}=\frac{r}{r+1}.$$ In fact the inequality is strict because otherwise $\rk K$ will be a multiple of $r+1$ which is a contradiction. so we have $$\slope{K}\geq\frac{r}{r+1}>\frac{r-2}{r-1}=\slope{E}.$$ Since $E$ is stable, there's no non-trivial map from $K$ to $E$. So the composition $K\longrightarrow V\longrightarrow G$ is non-trivial. Since $G$ is $\mu$-stable, we have $$\slope{K}\leq\slope{G}=1.$$ So now we have $$\frac{r}{r+1}<\slope{K}\leq 1.$$ 

If $$\frac{r}{r+1}<\slope{K}<1,$$ then $$0<\frac{\rk{K}-\deg K}{\rk{K}}<\frac{1}{r+1}.$$ Therefore we must have $\rk{K}>r+1$, which is absurd. Therefore, we must have $\slope{K}=1$. Now we consider the non-trivial map $\varphi: K\longrightarrow G$. It can be factored as $K\twoheadrightarrow\im{\varphi}\hookrightarrow G$. Since both $K$ and $G$ are $\mu$-stable, we have $$1=\slope{K}\leq\slope{\im{\varphi}}\leq\slope{G}=1.$$ So $\slope{\im{\varphi}}=1$, which implies that $\rk{\im{\varphi}}=\rk{G}$ and $\deg\im{\varphi}=\deg G$. We denote $Q=\coker{\varphi}$, then the support of $Q$ is 0 dimensional. From the exact sequence $K\longrightarrow G\longrightarrow Q\longrightarrow 0$, we have the exact sequence $\ext^1(Q,E)\longrightarrow\ext^1(G,E)\longrightarrow\ext^1(K,E)$. Since $Q$ has $0$ dimensional support and $E$ is locally free, we have $\ext^1(Q,E)=0$. Let $e\in\ext^1(G,E)$ be the extension class of $V$. Then the image of $e$ in $\ext^1(K,E)$ is given by the following pull back diagram: 
\begin{displaymath}
\xymatrix{0 \ar[r] & E \ar[r] \ar@{=}[d] & W \ar[r] \ar[d] & K \ar[r] \ar[d] \ar@{_{(}->}[ld] & 0\\
0 \ar[r] & E \ar[r] & V \ar[r] & G \ar[r] & 0.}
\end{displaymath}
Note that there's an injection from $K$ to $V$, therefore, there's an embedding from $K$ to $W$ which makes the first row split. This tells us that the image of $e$ in $\ext^1(K,E)$ is 0. From the above discussion we know that $e=0$, which is a contradiction. 

So $V$ is $\mu$-stable. 
\end{proof}


By applying the same method, we can prove

\begin{lemma}\label{l3-2}
For any $\mu$-stable locally free sheaf $G$ in $\fM_{X,H}(2v_2)$, and any extension given by a 2-dimensional subspace of $\ext^1(G,E)$ $$0\longrightarrow E^{\oplus 2}\longrightarrow F\longrightarrow G\longrightarrow 0,$$ the sheaf $F$ is a $\mu$-semistable locally free sheaf in $\fM_{X,H}(2v_1)$.
\end{lemma}

\begin{proof}
We consider the push-out diagram
\begin{displaymath}
\xymatrix{
 & 0 \ar[d] & 0 \ar[d] & & \\
 & E \ar@{=}[r] \ar[d] & E \ar[d] & & \\
0 \ar[r] & E^{\oplus 2} \ar[r] \ar[d] & F \ar[r] \ar[d] & G \ar[r] \ar@{=}[d] & 0\\
0 \ar[r] & E\oplus 0 \ar[r] \ar[d] & V \ar[r] \ar[d] & G \ar[r] & 0\\
 & 0 & 0 & & 
}
\end{displaymath}
Due to the fact that $F$ corresponds to a two dimensional subspace of $\ext^1(G,E)$, the extension $0\longrightarrow E\longrightarrow F\longrightarrow V\longrightarrow 0$ is a non-trivial extension. We want to show that $F$ is $\mu$-semistable. If not, let $K$ be a destabilizing sheaf of $F$. Then we have $$\slope{K}>\slope{F}=\frac{r-1}{r}.$$ Since $$\slope{E}=\frac{r-2}{r-1}<\slope{K},$$ we know there's no non-trivial map from $K$ to $E$. Therefore the composition $\varphi: K\longrightarrow F\longrightarrow V$ is a non-trivial map. It factor through $K\twoheadrightarrow\im\varphi\hookrightarrow V$. By the stability of $K$ and $V$, we have $$\slope{K}\leq\slope{\im\varphi}\leq\slope{V}.$$ So $$\frac{r-1}{r}<\slope{\im\varphi}\leq\frac{r}{r+1}.$$ 

If $$\frac{r-1}{r}<\slope{\im\varphi}<\frac{r}{r+1},$$ then we will have $$r<\frac{\rk\im\varphi}{\rk\im\varphi-\deg\im\varphi}<r+1,$$ which implies $\rk\im\varphi>2r$. Contradiction! 

If $\slope{\im\varphi}=\frac{r}{r+1}$, since $V$ is $\mu$-stable, we must have that $\rk\im\varphi=\rk V$ and $\deg\im\varphi=\deg V$. Let $Q=\coker\varphi$, then $Q$ is supported on a dimension 0 locus. The exact sequence $K\longrightarrow V\longrightarrow Q\longrightarrow 0$ gives an exact sequence $\ext^1(Q,E)\longrightarrow\ext^1(V,E)\longrightarrow\ext^1(K,E)$. Since the support of $Q$ has codimension 2, we have $\ext^1(Q,E)=0$. Let $e\in\ext^1(V,E)$ be the extension class corresponding to the non-trivial extension $0\longrightarrow E\longrightarrow F\longrightarrow V\longrightarrow 0$, then the image of $e$ in $\ext^1(K,E)$ is given by the pull back diagram
\begin{displaymath}
\xymatrix{
0 \ar[r] & E \ar[r] \ar@{=}[d] & H \ar[r] \ar[d] & K \ar[r] \ar[d] \ar@{_{(}->}[ld] & 0\\
0 \ar[r] & E \ar[r] & F \ar[r] & V \ar[r] & 0.
}
\end{displaymath}
From the construction of pull-back and the fact that $K$ is a subsheaf of $F$, we can see that the first row splits, which implies the image of $e$ in $\ext^1(K,E)$ is 0. So $e=0$. Contradiction!

This shows that $F$ must be $\mu$-semistable. 
\end{proof}


Now we want to prove the extension in the above lemma is essentially unique.

\begin{lemma}\label{3rdlemma}
For any $\mu$-stable locally free sheaf $G$ in $\fM_H(2v_2)$, $\dim\ext^1(G,E)=2$. In other words, there's only one such sheaf $F$ up to isomorphisms as in lemma \ref{l3-2}.
\end{lemma}
\begin{proof}
We have $\chi (G, E)=-\langle v(G), v(E)\rangle=-\langle 2v_2, v_0\rangle=-2$. However $\hom (G, E)=0$ since $\slope{G}>\slope{E}$. By Serre duality, we have $\ext^2(G, E)=\hom (E, G)$. 

Assume we have a nontrivial map $\varphi: E\longrightarrow G$. This map factor through $E\twoheadrightarrow\im{\varphi}\hookrightarrow G$. Since both $E$ and $G$ are $\mu$-stable, we have $$\slope{E}\leq\slope{\im{\varphi}}\leq\slope{G},$$ that is $$\frac{r-1}{r-2}\leq\slope{\im{\varphi}}\leq 1.$$ If both inequalities are strict, then $\rk\im{\varphi}\geq r-2$ which is a contradiction. Therefore we must have $\slope{\im{\varphi}}=1$ or $\slope{\im{\varphi}}=\frac{r-1}{r-2}$. 

If $\slope{\im{\varphi}}=1$, since $G$ is $\mu$-stable, the only possibility is $\rk{\im{\varphi}}=\rk{G}=2$ and $\deg\im{\varphi}=\deg G=2$. If we denote $Q=\coker{\varphi}$, then the support of $Q$ is 0 dimensional. From the exact sequence $E\longrightarrow G\longrightarrow Q\longrightarrow 0$, we have another exact sequence $\ext^1(Q,E)\longrightarrow\ext^1(G,E)\longrightarrow\ext^1(E,E)$. Since $Q$ has 0 dimensional support and $E$ is locally free, we have $\ext^1(Q,E)=0$. Due to the fact that $E$ is exceptional, we know that $\ext^1(E,E)=0$. So $\ext^1(G,E)=0$. Contradiction! 

If $\slope{\im{\varphi}}=\frac{r-1}{r-2}$, then $\rk{\im{\varphi}}$ must be an integer multiple of $r-2$. Since $\im{\varphi}$ is a quotient sheaf of $E$ we must have $\rk{\im{\varphi}}=r-1$ and $\deg \im{\varphi}=r-2$. Therefore $E=\im{\varphi}$. So we conclude that $\varphi:E\longrightarrow G$ must be an injection. However, when $r\geq 4$, we know that $\rk{E}=r-1>2=\rk{G}$. Contradiction! So for every $\mu$-stable locally free sheaf $G$, $\hom(E,G)=0$. Hence $\dim\ext^1(G,E)=2$ and we are done. 
\end{proof}


Now we are ready to finish the proof of proposition \ref{rk3}. 

\begin{proof}[Proof of proposition \ref{rk3}]

For every locally free $\mu$-stable sheaf $G$ with Mukai vector $2v_2$, by lemma \ref{l3-2} and \ref{3rdlemma}, we have associated it a locally free $\mu$-semistable sheaf $F$ with Mukai vector $2v_1$, which establishes a morphism from the locally free $\mu$-stable locus $U_2$ in $\fM_{X,H}(2v_2)$ to the moduli space $\fM_{X,H}^{\mu\textrm{-ss}}(2v_1)$ of $\mu$-semistable sheaves with Mukai vector $2v_1$. In fact, the extension sequence in lemma \ref{l3-2} holds in flat families. Let $V$ be any affine open subscheme of $U_2$. Let $\mathcal{G}$ be the universal sheaf on $V\times X$, and let $\mathcal{E}$ be the pullback of the exceptional bundle $E$ along the projection from $V\times X$ to $X$. Then we have the extension sequence $$0 \longrightarrow \mathcal{E} \longrightarrow \mathcal{F} \longrightarrow \mathcal{G} \longrightarrow 0$$ on $V\times X$, which induces a morphism from $V$ to $\fM_{X,H}^{\mu\textrm{-ss}}(2v_1)$. Furthermore, the restriction of the above extension sequence to every closed point in $\mathcal{F}$ is exactly the extension sequence in lemma \ref{l3-2}. Therefore, the morphism from all the affine subschemes of $U_2$ glue together to give a morphism $$\eta: U_2 \longrightarrow \fM_{X,H}^{\mu\textrm{-ss}}(2v_1).$$

Now we show this map is injective on closed points. From the exact sequence $$0\longrightarrow E^{\oplus 2}\longrightarrow F\longrightarrow G\longrightarrow 0,$$ we have $$0\longrightarrow\hom(E,E^{\oplus 2})\longrightarrow\hom(E,F)\longrightarrow\hom(E,G).$$
In lemma \ref{3rdlemma} we have proved $\hom(E,G)=0$, so $\dim\hom(E,F)=\dim\hom(E,E^{\oplus 2})=2$, which implies that there is only one way to get a quotient sheaf $G$. 

It's easy to see that the locus $U_1$ of $\mu$-stable locally free sheaves with Mukai vector $2v_1$ is also an open subscheme of $\fM_{X,H}^{\mu\textrm{-ss}}(2v_1)$. Therefore $\eta^{-1}(U_1)\cup U_2$ is an open subscheme of $U_2$, and hence an open subscheme of $\fM_{X,H}(2v_2)$. Since both $U_1$ and $U_2$ are smooth, the restriction of $\eta$ on $\eta^{-1}(U_1)\cap U_2$ is an isomorphism onto its image, which identifies two smooth open subschemes of $\fM_{X,H}(2v_1)$ and $\fM_{X,H}(2v_2)$. So we can conclude that $\fM_{X,H}(2v_1)$ and $\fM_{X,H}(2v_2)$ are birational. 
\end{proof}


Finally, we deal with two exceptional cases $r=2$ and $r=3$. We have two similar birational results as in the above general case. 

\begin{proposition}\label{exp1}
For any given $s\geq 1$, let $X$ be a projective K3 surface with $\pic{X}=\mathbb{Z}[H]$ where $H$ is an ample line bundle with $H^2=50s-28$. Consider Mukai vectors
\begin{eqnarray*}
v_0&=&(5, 2H, 20s-11),\\
v_1&=&(7, 3H, 32s-18),\\
v_2&=&(2, H, 12s-17).
\end{eqnarray*}
Then the moduli spaces of semi-stable sheaves $\fM_{X,H}(2v_1)$ is birational equivalent to $\fM_{X,H}(2v_2)$. \qed
\end{proposition}

\begin{proposition}\label{exp2}
For any given $s\geq 1$, let $X$ be a projective K3 surface with $\pic{X}=\mathbb{Z}[H]$ where $H$ is an ample line bundle with $H^2=98s-60$. Consider Mukai vectors
\begin{eqnarray*}
v_0&=&(7, 2H, 28s-17),\\
v_1&=&(10, 3H, 44s-27),\\
v_2&=&(3, H, 16s-10).
\end{eqnarray*}
Then the moduli spaces of semi-stable sheaves $\fM_{X,H}(2v_1)$ is birational equivalent to $\fM_{X,H}(2v_2)$. \qed
\end{proposition}

The proof of the above two propositions are completely parallel to the proof of proposition \ref{rk3}. We will not repeat them here. Combining corollary \ref{_def_cor_}, proposition \ref{rk3}, \ref{exp1}, \ref{exp2}, we have proved theorem \ref{_thm_main_}.

\section{Rank 0 Case} 

In this section we will prove the following birational equivalence between a moduli space of rank 0 sheaves and a moduli space of rank 2 sheaves. 

\begin{theorem}\label{rank0}
Let $X$ be a K3 surface and $H$ be an ample line bundle. Assume $\pic{X}=\mathbb{Z}[H]$. Let $v_1=(0, 2H, -2)$ and $v_2=(2, 0, -H^2)$. Then the moduli spaces of sheaves $\fM_{X,H}(v_1)$ and $\fM_{X,H}(v_2)$ are birational. 
\end{theorem}

\begin{proof}
Let $v_3=v_2\cdot\ch{H}=(2, 2H, 0)$. Note that by tensoring the line bundle $H$, we get an isomorphism between $\fM_{X,H}(v_2)$ and $\fM_{X,H}(v_3)$. Therefore we only need to prove $\fM_{X,H}(v_1)$ is birational to $\fM_{X,H}(v_3)$. 

Let $U=\{\ i_*Q\in\fM_{X,H}(v_1)\ |\ Q\ \mbox{is a line }\mbox{bundle }\mbox{supported }\mbox{on a }\mbox{smooth }\mbox{curve}\ C \in |2H|\}$. Then $U$ is an open subscheme of $\fM_{X,H}(v_1)$. We know that $g(C)=1+\frac{1}{2}C^2=1+2H^2$. Furthermore, by Grothendieck-Riemann-Roch formula, we know that $$i_*(\ch{Q}\cdot\td{C})=\ch{i_*Q}\cdot\td{X},$$ that is, $$i_*((1, \deg Q)\cdot(1, -2H^2))=(0, 2H, -2)\cdot(1, 0, 2).$$ We can conclude that $\deg Q=2H^2-2$. 

For every $i_*Q\in U\subset\fM_{X,H}(v_1)$, we construct a sheaf $F$ as a non-trivial extension $$0\longrightarrow\mathcal{O}^{\oplus 2}\longrightarrow F\longrightarrow i_*Q\longrightarrow 0.$$

First of all, we have to show that, for a generic $i_*Q\in U$, there's only one way to build up such an extension $F$. In order to prove this, we need to show that a generic $i_*Q\in U$ satisfies $\dim\ext^1(i_*Q, \mathcal{O})=2$. 

By Serre duality, we know that 
\begin{eqnarray*}
\ext^1(i_*Q, \mathcal{O})&=&\ext^1(\mathcal{O}, i_*Q\otimes K_X)^{\vee}=H^1(X, i_*Q\otimes K_X)^{\vee}\\
&=&H^1(C, Q\otimes i^*K_X)^{\vee}=H^0(C, Q^{\vee}\otimes N_{C/X}).
\end{eqnarray*}
We can find that $\deg(Q^{\vee}\otimes N_{C/X})=\deg(Q^{\vee}\otimes K_C)=2H^2+2$. By Riemann-Roch formula, we know that $$\chi(Q^{\vee}\otimes N_{C/X})=1-(2H^2+1)+(2H^2+2)=2.$$

We need to prove that for a generic $i_*Q\in U$, $\dim H^0(C, Q^{\vee}\otimes N_{C/X})=2$. By upper-semi-continuity theorem, we only need to show that there exists at least one $i_*Q\in U$, which makes $\dim H^0(C, Q^{\vee}\otimes N_{C/X})=2$. 

Since $Q$ is an arbitrary line bundle on an arbitrary smooth curve $C\in|2H|$, we only need to prove the following lemma:

\begin{lemma}\label{open1}
Let $C$ be a smooth curve and $L$ be a line bundle on $C$ with $\chi(L)=2$. If $L$ is general, then $\dim H^0(L)=2$ and $\dim H^1(L)=0$. 
\end{lemma}

\begin{proof}
In contrary we assume the infimum of $\dim H^0(L)$ of such $L$'s is at least 3. Pick any $p\in C$ which is not a base point of $L$. From the exact sequence $$0\longrightarrow L(-p)\longrightarrow L\longrightarrow \mathcal{O}_p\longrightarrow 0$$ and its associated long exact sequence 
\begin{eqnarray*}
0&\longrightarrow& H^0(L(-p))\longrightarrow H^0(L)\longrightarrow \mathbb{C}_p\\
&\longrightarrow& H^1(L(-p))\longrightarrow H^1(L)\longrightarrow 0
\end{eqnarray*}
we know that $\dim H^0(L(-p))=\dim H^0(L)-1$, and $\dim H^1(L(-p))=\dim H^1(L)\geq 1$. 

By Serre duality, we have $H^1(L(-p))=H^0(K_C\otimes L^{\vee}(p))^{\vee}$. We can again pick a point $q\in C$ which is not a base point of $K_C\otimes L^{\vee}(p)$. Then from the short exact sequence $$0\longrightarrow K_C\otimes L^{\vee}(p-q)\longrightarrow K_C\otimes L^{\vee}(p)\longrightarrow\mathcal{O}_q\longrightarrow 0$$ and its associated long exact sequence 
\begin{eqnarray*}
0&\longrightarrow& H^0(K_C\otimes L^{\vee}(p-q))\longrightarrow H^0(K_C\otimes L^{\vee}(p))\longrightarrow \mathbb{C}_q\\
&\longrightarrow& H^1(K_C\otimes L^{\vee}(p-q))\longrightarrow H^1(K_C\otimes L^{\vee}(p))\longrightarrow 0
\end{eqnarray*}
we know that $$\dim H^0(K_C\otimes L^{\vee}(p-q))=\dim H^0(K_C\otimes L^{\vee}(p))-1$$ and $$\dim H^1(K_C\otimes L^{\vee}(p-q))=\dim H^1(K_C\otimes L^{\vee}(p)).$$ Again by Serre duality, we know that $$\dim H^1(L(-p+q))=\dim H^1(L(-p))-1=\dim H^1(L)-1$$ and $$\dim H^0(L(-p+q))=\dim H^0(L(-p))=\dim H^0(L)-1.$$

So we see that the existence of the line bundle $L(-p+q)$ conflicts the assumption. Therefore, there must be a line bundle $L$ with $\chi(L)=2$, such that $\dim H^0(L)=2$. By the upper-semi-continuity theorem, we know that this is true for a generic line bundle on any smooth curve $C$. 
\end{proof}

Back to the proof of theorem \ref{rank0}. We know for a generic $i_*Q\in U$, we have $\dim\ext^1(i_*Q, \mathcal{O})=2$. Therefore, there's only one way to produce a sheaf $F\in\fM_{X,H}(v_3)$ from $i_*Q$ via the extension $$0\longrightarrow\mathcal{O}^{\oplus 2}\longrightarrow F\longrightarrow i_*Q\longrightarrow 0.$$

Now we have to prove $F$ is torsion free. In fact, we want to prove, if $i_*Q\in U$ is generic, then $F$ is locally free. We need the following two lemmas:

\begin{lemma}
Under the above conditions, the sheaf $F$ is locally free if and only if $Q^{\vee}\otimes N_{C/X}$ is base point free. 
\end{lemma}
\begin{proof}
Note that $\chom(i_*Q, \mathcal{O}^{\oplus 2})=0$. By local to global spectral sequence, we have an isomorphism $$\ext^1(i_*Q, \mathcal{O}^{\oplus 2})\cong H^0(\cext^1(i_*Q, \mathcal{O}^{\oplus 2}))=H^0(\cext^1(i_*Q, \mathcal{O})\oplus\cext^1(i_*Q, \mathcal{O})).$$ Therefore, to prove the sheaf $F$ is locally free, we only need to show the restriction of the two global sections of $\cext^1(i_*Q, \mathcal{O})$ to every local ring, provide a free local extension of the two sheaves. 

Let $p$ be a closed point of $C$. Over $\textrm{Spec}\,\mathcal{O}_p$, we can assume the curve is cut by one equation $g\in\mathcal{O}_p$, and resolve the torsion sheaf $i_*Q$ as $$0\longrightarrow\mathcal{O}_p\stackrel{\cdot\,g}{\longrightarrow}\mathcal{O}_p \longrightarrow\mathcal{O}_p/(g)\longrightarrow 0.$$
Apply the functor $\hom(-, \mathcal{O}_p^{\oplus 2})$ we get the long exact sequence $$0\longrightarrow\hom(\mathcal{O}_p, \mathcal{O}_p^{\oplus 2})\longrightarrow\hom(\mathcal{O}_p, \mathcal{O}_p^{\oplus 2})\longrightarrow\ext^1(\mathcal{O}_p/(g), \mathcal{O}_p^{\oplus 2})\longrightarrow 0.$$
Therefore, let $(\varphi_1, \varphi_2)$ be the restriction of two linearly independent sections of $\cext^1(i_*Q, \mathcal{O})$ to the local ring $\mathcal{O}_p$, then
$(\varphi_1, \varphi_2)\in\cext^1(i_*Q, \mathcal{O}^{\oplus 2})_p=\ext^1(\mathcal{O}_p/(g), \mathcal{O}_p^{\oplus 2})$. They are images of functions $(f_1, f_2)\in\hom(\mathcal{O}_p, \mathcal{O}_p^{\oplus 2})$. Then the stalk of the sheaf $F$ over $\textrm{Spec}\,\mathcal{O}_p$ is given by the push out diagram
\begin{displaymath}
\xymatrix{
0 \ar[r] & \mathcal{O}_p \ar[d]_g \ar[r]^{(f_1, f_2)} & \mathcal{O}_p^{\oplus 2} \ar[r] \ar[d] & \mathcal{O}_p/(g) \ar[r] \ar@{=}[d] & 0\\
0 \ar[r] & \mathcal{O}_p \ar[r] & F_p \ar[r] & \mathcal{O}_p/(g) \ar[r] & 0
}
\end{displaymath}
Therefore, the stalk $F_p=\mathcal{O}_p^{\oplus 3}/(f_1, f_2, g)$. Since $g$ is the defining equation of the curve $C$, $g(p)=0$. So $F_p$ is a free $\mathcal{O}_p$ module if and only if $f_1(p)$ and $f_2(p)$ are not simultaneously 0. 

Note that $(\varphi_1, \varphi_2)$ is the restriction of $(f_1, f_2)$ on the curve $C$. So $F_p$ is free when $\varphi_1(p)$ and $\varphi_2(p)$ are not simultaneously 0. Therefore, as long as $\varphi_1$ and $\varphi_2$ don't have common zeroes along $C$, the sheaf $F$ is locally free. That is to say, we want the sheaf $\cext^1(i_*Q, \mathcal{O})$ to be base point free for a generic $i_*Q$. 

Furthermore, we can also observe from the above that the sheaf $\cext^1(i_*Q, \mathcal{O})$ is canonically isomorphic to $i_*\chom(Q\otimes\mathcal{O}_C(-C), \mathcal{O}_C)=i_*(Q^{\vee}\otimes N_{C/X})$. So we have $H^0(X, \cext^1(i_*Q, \mathcal{O}))=H^0(C, Q^{\vee}\otimes N_{C/X})$. Therefore, the extension $F$ is locally free if and only if $Q^{\vee}\otimes N_{C/X}$ is base point free. 
\end{proof}

We have already proved that for a generic $i_*Q\in U$, $Q^{\vee}\otimes N_{C/X}$ has two dimensional global sections. The next step is to prove that under this condition, a generic $Q^{\vee}\otimes N_{C/X}$ is base point free.

\begin{lemma}\label{open2}
Let $C$ be a smooth curve and $L$ be a line bundle on $C$ with $\chi(L)=2$. For a general $L$, it is base point free. 
\end{lemma}
\begin{proof}
Without loss of generality, we only need to consider all line bundles $L$ with $\dim H^0(L)=2$ and $\dim H^1(L)=0$. Note that base point freeness is an open condition in a flat family of line bundles with constant dimensional cohomology groups. Therefore we only need to show that there is at least one such line bundle $L$ which is base point free. 

First of all we can choose a line bundle $M$ such that $\dim H^0(M)=1$ and $\dim H^1(M)=0$. In fact, for any line bundle $L$ with $\dim H^0(L)=2$ and $\dim H^1(L)=0$, let $p$ be a point not in the base locus of $L$, then $M=L(-p)$ does the job. We denote the generator of $H^0(M)$ by $\sigma$. 

We claim that there must be a point $q$ which makes $M(q)$ base point free. In fact, from the exact sequence of cohomology groups, we know $\dim H^0(M(q))=2$ and $\dim H^1(M(q))=0$ for any $q\in C$. We assume the contrary that for every $q\in C$, $M(q)$ has a base point. Then the base point must be in the zero locus of $\sigma$. Since the zero set of the section $\sigma$ is finite, we conclude that there must be a point $r\in\sigma^{-1}(0)$, which is the common base point of $M(q)$'s for all $q\in C$. In particular, $r$ is the base point of $M(r)$ which is absurd. This proves the lemma. 
\end{proof}

From the above two lemmas, we see that for a generic $i_*Q\in U$, the extension sheaf $F$ we constructed is a locally free sheaf. 

Now let's prove the sheaf $F$ is stable. 

Assume that $F$ can be destabilized by a subsheaf $L$. Since $F$ is locally free of rank 2, without loss of generality we can further assume $\rk{L}=1$ and $L$ is locally free. From $\slope{L}\geq\slope{F}$, we have $c_1(L)\geq H$. We assume that $L=\mathcal{O}(H)$, since otherwise $\mathcal{O}(H)$ also destabilizes $F$. 

Obviously $\hom (L, \mathcal{O}^{\oplus 2})=0$ since $\slope{L}>0$. We also have $$\hom(L, i_*Q)=H^0(X, L^{\vee}\otimes i_*Q)=H^0(C, i^*L^{\vee}\otimes Q).$$ Since $\deg(i^*L^{\vee}\otimes Q)=\deg L^{\vee}+\deg Q=-L\cdot C+\deg Q=-2H^2+(2H^2-2)=-2$, we know that $\hom(L, i_*Q)=H^0(C, i^*L^{\vee}\otimes Q)=0$. Therefore $\hom(L, F)=0$. Contradiction! 

Finally, let's prove for every such sheaf $F$ that we obtained by extension $$0\longrightarrow\mathcal{O}^{\oplus 2}\longrightarrow F\longrightarrow i_*Q\longrightarrow 0,$$ there's only one sheaf $i_*Q\in U$ with $\dim H^0(C, Q^{\vee}\otimes N_{C/X})=2$ which produces $F$. In fact, in the long exact sequence of cohomology groups $$0\longrightarrow H^0(X, \mathcal{O}^{\oplus 2})\longrightarrow H^0(X, F)\longrightarrow H^0(X, i_*Q),$$ we already know that $$H^0(X, i_*Q)=H^0(C, Q)=H^1(C, Q^{\vee}\otimes K_C)^{\vee}=H^1(C, Q^{\vee}\otimes N_{C/X})^{\vee}=0.$$ So $\dim H^0(F)=\dim H^0(\mathcal{O}^{\oplus 2})=2$. Therefore $i_*Q$ is also uniquely determined by $F$. 

In the same way as in the proof of proposition \ref{rk3}, we can see that the above extension can also be done in flat falilies, therefore gives a morphism $$\eta: V\longrightarrow \fM_{X,H}(v_3)$$ where  $V$ is the open subscheme of $\fM_{X,H}(v_1)$, which parametrizes all line bundles supported on smooth curves in the linear system $|2H|$ satisfying the conditions of lemmas \ref{open1} and \ref{open2}. And the image of $\eta$ lies in the smooth open subscheme of $\fM_{X,H}(v_3)$ which parametrizes locally free stable sheaves with Mukai vector $v_3$. The above discussion also shows that $\eta$ identifies $V$ with its image. Since the moduli spaces $\fM_{X,H}(v_1)$ and $\fM_{X,H}(v_3)$ are both irreducible and of the same dimension, $\eta$ establishes a birational morphism from $\fM_{X,H}(v_1)$ to $\fM_{X,H}(v_3)$. By tensoring with the line bundle $H^\vee$ we see that $\fM_{X,H}(v_1)$ and $\fM_{X,H}(v_2)$ are birational.

\end{proof}

\bibliographystyle{alpha}
\bibliography{references}

\end{document}